\documentclass[12pt,french]{article}

\usepackage[applemac]{inputenc}
\usepackage{ae, aeguill}
\usepackage[francais]{babel} 
\usepackage{amssymb}

\def \N{{\mathbb N}}

\def \Z{{\mathbb Z}}

\def \P{{\mathbb P}}
\def \E{{\mathbb E}}

\def\1{1\!\! 1}

\def\±{\displaystyle}

\def\vs2{\vspace{2mm}}
\def\vs3{\vspace{3mm}}
\def\vs5{\vspace{5mm}}

\def \norm#1#2{\| #1 \|_{#2}}

\begin{document}

\vspace{1cm}
\centerline{\bf A simple proof of a recurrence theorem for random walks in $\Z^{2}$}

\vspace{2mm}
\centerline{Jean-Marc Derrien}

\vspace{2mm}
\centerline{\it October, 2006}

\vspace{1cm}
\noindent{\bf Abstract }{\it In this note, we prove without using Fourier analysis
that the symmetric square integrable
random walks in $\Z^{2}$ are recurrent.}

\vspace{1cm}
George Pòlya related in [3, pp. 582-583] an incident that
enables him to formulate the question 
of recurrence for random walks: during a stroll through the woods,
he felt embarrassed because he met ``certainly much too often'' a student with his
girlfriend.

\vspace{3mm}
Let $x$ be an element of $\Z^{2}$.
A random walk in $\Z^{2}$ starting at $x$
is a sequence $(S_{n})_{n\in\N}$ of random variables
such that
$$
S_{0}=x
\quad{\rm and}\quad
S_{n}=x+X_{1}+X_{2}+\cdots+X_{n}\; , \quad n\in \N\; ,
$$
where $(X_{n})_{n\in \N}$ is a sequence 
of independent and identically distributed (i.i.d.), $\Z^{2}$-valued random variables.
$(S_{n})_{n\in \N}$ is characterized
in law by  $x$ and the law of $X_{1}$.

\vspace{3mm}
In this note, we give an elementary proof of the following result.

\vspace{2mm}
\noindent{\bf Theorem}
{\it If $(S_{n}^{(1)})_{n\in \N}$ and $(S_{n}^{(2)})_{n\in \N}$
are two square integrable, independent and identically distributed random walks
in $\Z^{2}$,
then
$$
S_{n}^{(1)}=S_{n}^{(2)}
\quad
\mbox{infinitely often}
$$
with probability one.}

\vspace{3mm}
This theorem is a straightforward consequence of the following proposition
(first proved in [2] in a more general context)
applied to
$(S_{n}:=S_{n}^{(1)}-S_{n}^{(2)})_{n\in \N}$.

\vspace{2mm}
\noindent{\bf Proposition}
{\it If $(S_{n})_{n\in \N}$ is a symmetric, square integrable random walk
starting at 0 in $\Z^{2}$,
then
$$
S_{n}=0
\quad
\mbox{infinitely often}
$$
with probability one.}

\vspace{2mm}
\noindent{\it Proof }
It is classical that,
in order to prove this proposition,
we have only to establish that
$$
\sum_{n=0}^{+\infty} \P [S_{n}=0]=+\infty
$$
(see, for instance, [1]).

\vspace{2mm}
One can write
$$
S_{n}=X_{1}+X_{2}+\cdots+X_{n} \; ,
$$
where $(X_{k})_{k\geq 1}$ is a sequence of
i.i.d., square integrable, $\Z^{2}$-valued  random variables
with
$X_{1}\equiv -X_{1}$
in distribution.

\vspace{2mm}
Since the square integrable random walk $(S_{n})_{n\in\N}$
is symmetric, it is centered and we have
$$
\E ( \norm{S_{n}}{2}^{2})
=
\sum_{k=1}^{n} \E (\norm{X_{k}}{2}^{2})
+
\sum_{1 \leq k < l \leq n} \E (X_{k})\cdot \E (X_{l})
=
n\; \E ( \norm{X_{1}}{2}^{2})
$$
($||\cdot||_{2}$ denotes the euclidean norm).

\vspace{2mm}
The symmetry also gives
\begin{eqnarray*}
\P [S_{2n} =0 ]
&=&
\sum_{x\in \Z^2} \P [S_{2n}=0\; |\ S_{n}=x] \; \P [S_{n}=x]
\\
&=&
\sum_{x\in \Z^{2}} \P [-X_{n+1}-X_{n+2}-\cdots -X_{2n}=x]\; 
\P [S_{n}=x]
=\sum_{x\in \Z^2} \P [S_{n}=x]^2\; .
\end{eqnarray*}

\vspace{2mm}
Hence if we introduce
$$
B_{n}:=
\left\{x\in \Z^{2}\; :\ \norm{x}{2}^{2}< 2n\;
\E \left( \norm{X_{1}}{2}^{2}\right) \right\}\; ,
$$
we deduce from Cauchy-Schwarz's and
Markov's inequalities that, if $n$ is large enough,
\begin{eqnarray*}
\P [S_{2n} =0 ]
&\geq&
{1\over{
| B_{n} |
}}
\left(\sum_{x\in B_{n}}
\P [S_{n}=x] \right)^{2}
\\
&\geq&
{C \over
{n}
}
\bigg( 1-\P [ \norm{S_{n}}{2}^{2} \geq 2 n\; \E (\norm{X_{1}}{2}^{2})] \bigg)^{2}
\\
&\geq&
{C\over{n}}
\left(
1- {{\E \left( \norm{S_{n}}{2}^{2} \right)}
\over{2 n\; \E (\norm{X_{1}}{2}^{2})}}
\right)^{2}
=
{C \over {4n}}
\; ,
\end{eqnarray*}
where $C>0$ depends only on $\E (\norm{X_{1}}{2}^{2})$.
The proposition follows.

\hfill$\bullet$

 \vspace{1cm}
 \noindent{\bf Acknowledgement }I am grateful to Emmanuel Lesigne
 for useful remarks.

\vspace{1cm}
\noindent{REFERENCES}

\vspace{2mm}
\noindent [1] L. Breiman, {\it Probability},
Addison-Wesley, Reading, 1968.

\vspace{1mm}
\noindent [2] K. L. Chung and W. H. J. Fuchs, On the distribution
of values of sums of random variables, {\it Mem. Amer.  Math. Soc.} {\bf 6}, 1951. 

\vspace{1mm}
\noindent [3] G. Pòlya, {\it Collected papers}, {\bf volume IV}, The MIT
Press, Cambridge, Massachusetts, 1984.

\vspace{1cm}
\noindent {\it Département de Mathématiques. Université de Brest. France.}

\noindent{\it Jean-Marc.Derrien@univ-brest.fr}
\end{document}